\documentclass[11pt]{amsart}
\usepackage{amssymb,pstricks,pstcol,pst-plot}

\definecolor{Red}{cmyk}{0,1,1,0}
\definecolor{OrangeRed}{cmyk}{0,0.6,1,0}            
\definecolor{DarkBlue}{cmyk}{1,1,0,0.20}
\definecolor{Black}{cmyk}{0,0,0,1}
\definecolor{Violet}{cmyk}{0.79,0.88,0,0}
\definecolor{Purple}{cmyk}{0.45,0.86,0,0}
\definecolor{myblue}{rgb}{0.66,0.78,1.00}

\parskip=\smallskipamount
\numberwithin{equation}{section}

\newtheorem{theorem}{Theorem}[section]
\newtheorem{lemma}[theorem]{Lemma}
\newtheorem{corollary}[theorem]{Corollary}
\newtheorem{proposition}[theorem]{Proposition}

\theoremstyle{definition}
\newtheorem{definition}[theorem]{Definition}
\newtheorem{example}[theorem]{Example}

\newtheorem{problem}[theorem]{Problem}
\newtheorem{remark}[theorem]{Remark}

\newcommand{\B}{\mathbb{B}}
\newcommand{\C}{\mathbb{C}}
\newcommand{\D}{\mathbb{D}}

\newcommand{\N}{\mathbb{N}}

\newcommand{\Z}{\mathbb{Z}}

\newcommand{\R}{\mathbb{R}}
\newcommand{\T}{\mathbb{T}}

\newcommand{\cA}{\mathcal{A}}
\newcommand{\cB}{\mathcal{B}}

\newcommand{\cC}{\mathcal{C}}
\newcommand{\cD}{\mathcal{D}}

\newcommand{\cH}{\mathcal{H}}

\newcommand{\cO}{\mathcal{O}}

\newcommand{\cU}{\mathcal{U}}

\def\e{\epsilon}

\def\dim{{\rm dim}\,}

\def\di{\partial}
\def\dibar{\bar\partial}
\def\bs{\backslash}

\def\wt{\widetilde}

\begin{document}
\title[Manifolds of holomorphic mappings]
{Manifolds of holomorphic mappings from strongly pseudoconvex domains}
\author{Franc Forstneri\v c}
\address{Institute of Mathematics, Physics and Mechanics, 
University of Ljubljana, Jadranska 19, 1000 Ljubljana, Slovenia}
\email{franc.forstneric@fmf.uni-lj.si}
\thanks{Research supported by grants P1-0291 and J1-6173, Republic of Slovenia.}

%
%    General info
%

\date{November 26, 2006} 

\begin{abstract}
Let $D$ be a bounded strongly pseudoconvex domain in a Stein manifold, 
and let $Y$ be a complex manifold. We show that many classical spaces 
of maps $\bar D \to Y$ which are holomorphic in $D$ are 
infinite dimensional complex manifolds which are modeled on 
locally convex topological vector spaces (Banach, Hilbert or Fr\'echet).  
This holds in particular for H\"older and Sobolev spaces of holomorphic maps.
\end{abstract}
\maketitle

%
%
%  INTRODUCTION
%
%
\section{Introduction}
Given complex manifolds (or complex spaces) $X$ and $Y$, it is
a natural question whether the space $\cH(X,Y)$ of all holomorphic mappings 
$X\to Y$ is also a complex manifold (resp.\ a complex space). 
If $X$ is a compact complex space without boundary then $\cH(X,Y)$ 
is a finite dimensional complex space which can be identified with 
an open subset in the Douady space 
$\cD(X\times Y)$, \cite[Theorem 1.5]{Campana}. 
The set of all holomorphic maps from a noncompact manifold in 
general does not admit any particularly nice structure.

In \S 2 of this paper we prove the following results.

\begin{theorem}
\label{manifolds}
Let $D$ be a relatively compact strongly pseudoconvex domain in a 
Stein manifold, and let $Y$ be a complex manifold. 
\begin{itemize}
\item[(i)] The H\"older space $\cA^{k,\alpha}(D,Y)=\cC^{k,\alpha}(\bar D,Y) \cap \cH(D,Y)$ 
is a complex Banach manifold for every $k\in\Z_+$ and $0\le \alpha <1$.
\item[(ii)] $\cA^\infty(D,Y)= \cC^\infty(\bar D,Y) \cap \cH(D,Y)$ 
is a complex Fr\'echet manifold. 
\item[(iii)] The Sobolev space $L^{k,p}_\cO(D,Y)= L^{k,p}(\bar D,Y) \cap \cH(D,Y)$ 
is a complex Banach manifold for $k\in\N$, $p\ge 1$ and $kp > \dim_{\! \R} D$ 
(resp.\ a complex Hilbert manifold if $p=2$).
\end{itemize}
If $L(D,Y)$ denotes any of the above manifolds of maps then 
the tangent space $T_f L(D,Y)$ at a point $f\in L(D,Y)$
is $L_h(D,f^*TY)$, the space of sections of class $L(D)$
of the complex vector bundle $h\colon f^*TY\to \bar D$.
If $D$ is contractible, or if $\dim D=1$, then 
$T_f L(D,Y) \approx L(D,\C^m)$ with $m=\dim Y$.
\end{theorem}

The analogous conclusions hold if $\bar D$ is a compact 
complex manifold with Stein interior $D$ and smooth strongly pseudoconvex 
boundary $bD$; according to  Heunemann \cite{Heunemann3} and Ohsawa \cite{Ohsawa}
(see also Catlin \cite{Catlin}) such $\bar D$ embeds as a 
smoothly bounded strongly pseudoconvex domain in a Stein manifold.

The special case of Theorem \ref{manifolds} (i) 
with $\alpha=0$ was proved recently in~\cite{BDF2}.

Before proceeding we recall some of the known results on this subject.
If $M$ is a compact smooth manifold, possibly with boundary, and $Y$ is 
a smooth manifold without boundary then the space of H\"older maps $\cC^{k,\alpha}(M,Y)$ 
$(k\in\Z_+,\ 0\le \alpha<1$), the space of Sobolev maps $L^{k,p}(M,Y)$ 
$(k\in \N,\ 1\le p<\infty,\ kp>\dim M)$, as well as many other natural mapping
spaces, are infinite dimensional Banach manifolds (Fr\'echet if we consider
$\cC^\infty$ maps). For results in this direction (also for 
sections of smooth fiber bundles over compact smooth manifolds) see 
Palais \cite{Palais1,Palais2}, El\u\i asson \cite{Eliasson}, 
Krikorian \cite{Krikorian}, Penot \cite{Penot}, Riddell \cite{Riddell};
for noncompact source manifolds see also Cantor \cite{Cantor}.
The cited works deal with important foundational questions of 
nonlinear global analysis and provide the basic framework for the
study of global partial differential operators and the calculus of variations.

These results do not apply directly in the holomorphic category,
the main reason being that it is much more difficult (and in general
impossible) to find holomorphic tubular neighborhoods. 
Besides the already mentioned work of Douady
and others on compact complex subspaces of a complex space, we 
are only aware of a few scattered results.
Ivashkovich and Shevchishin considered certain manifolds of 
pseudoholomorphic maps from Riemann surfaces to almost complex manifolds 
\cite{IS1,IS2}. Quite recently, Lempert \cite{Lempert} proved that for a compact 
smooth manifold $M$ and a complex manifold $Y$
the space $\cC^k(M,Y)$ is an infinite dimensional complex 
Banach manifold (Fr\'echet if $k=\infty$).  
It is expected that, like their real counterparts, 
these {\em generalized loop spaces} could be useful in the 
study of geometric properties of the target manifold. 

Our proof of Theorem \ref{manifolds} is based on the following result 
of possible independent interest,
concerning the existence of tubular Stein neighborhoods of holomorphic graphs 
with continuous boundary values over strongly pseudoconvex Stein domains.

\begin{theorem}
\label{Stein-nbds}
Assume that $h\colon X\to S$ is a holomorphic submersion of a complex manifold $X$
onto a Stein manifold $S$, $D\Subset S$ is a strongly pseudoconvex domain 
with $\cC^2$ boundary in $S$, and $f\colon \bar D\to X$ is a continuous section  
of $h$  which is holomorphic in $D$. 

There exists a holomorphic vector bundle
$\pi\colon E\to U$ over an open set $U\subset S$ containing $\bar D$,
and for every open set $\Omega_0\subset X$ containing $f(\bar D)$ there 
exist a Stein open set $\Omega$ in $X$ with $f(\bar D)\subset \Omega\subset \Omega_0$ 
and a biholomorphic map $\Theta \colon \Omega \to \Theta(\Omega) \subset E$ 
which maps the fiber $\Omega_z = h^{-1}(z) \cap \Omega$ over any point 
$z\in h(\Omega) \subset U$ biholomorphically onto an open convex 
set $\wt\Omega_z=\Theta(\Omega_z)$ in the fiber $E_z=\pi^{-1}(z)$.
\end{theorem}

If $\Theta$ is as in Theorem \ref{Stein-nbds}
then the map $g\to \Theta\circ g$ induces an isomorphism between the 
space of sections of the restricted submersion 
$X_{\bar D} = h^{-1}(\bar D) \to\bar D$ which are sufficiently uniformly 
close to $f$, and the space of sections of the complex vector bundle $E_{\bar D} \to \bar D$ 
which are close to $\Theta \circ f$. For each of the mapping spaces in 
Theorem \ref{manifolds} this provides a local holomorphic chart around an 
element $f$, and it is easily seen that such charts are holomorphically compatible 
and hence define a complex manifold structure on the relevant 
space of sections of $X_{\bar D} \to \bar D$. 
Further details of the proof of Theorem \ref{manifolds} are given in \S 2.

In the special case when the section $f$ in Theorem \ref{Stein-nbds} 
extends holomorphically to a neighborhood of $\bar D$, the conclusion 
of Theorem \ref{Stein-nbds} follows by combining Siu's theorem \cite{Siu} 
(or Schneider's theorem \cite{Sch}) with the Docquier-Grauert tubular neighborhood theorem
\cite{DG}; in this case $\Theta$ can be chosen to map $f$ onto the zero section of $E$.

The following is an immediate corollary to Theorem  \ref{Stein-nbds}.
Examples in \S 5 show that the conclusion fails in general for 
{\em images} of holomorphic maps.

\begin{corollary}
\label{cor1}
If $D\Subset S$ are as in Theorem \ref{Stein-nbds}, $Y$ is a
complex manifold and $f\colon \bar D\to Y$ is a continuous map
which is holomorphic in $D$ then the graph $G_f=\{(z,f(z))\colon z\in \bar D\}$ 
admits a basis of open Stein neighborhoods in $S\times Y$.
\end{corollary}

There are at least two possible ways to prove Theorem 
\ref{Stein-nbds}. One approach is by plurisubharmonic functions
and Grauert's characterization of Stein manifolds 
\cite[Theorem 2]{Grauert} (I wish to thank the referee for indicating this possibility).
To this end one needs
\begin{itemize}
\item[(i)] a strongly plurisubharmonic function $\rho$ in a neighborhood
of the compact set $f(\bar D)$ in $X$, and
\item[(ii)] a nonnegative (weakly) plurisubharmonic function $\tau$
in a neighborhood of $f(\bar D)$ which vanishes precisely on $f(\bar D)$.
\end{itemize}
For small enough $\e>0$ the function $\rho+\frac{1}{\e - \tau}$
is then a strongly plurisubharmonic exhaustion function on 
$\Omega_\e = \{\tau <\e\} \subset X$ and hence $\Omega_\e$ is Stein.
Functions $\rho$ and $\tau$ with these properties can be found
by classical methods.

We choose to present a different proof which is based on 
the method of {\em holomorphic sprays} developed in \cite{BDF1,BDF2}. 
This method explicitly produces a holomorphic vector bundle 
$\pi\colon E\to U$ as in Theorem \ref{manifolds} whose total space $E$ contains
a  biholomorphic copy of a Stein neighborhood of $f(\bar D)$.
We give a brief outline; for the details see Sec.\ \S 4.

In \cite{BDF2} the authors constructed a dominating spray
with a given central section, holomorphic over $D$ and
continuous over $\bar D$, by successively gluing sprays 
over small subsets of $\bar D$. 
The key ingredient of this construction is a 
Cartan type splitting lemma  with control up to
the boundary \cite[Theorem 3.2]{BDF1}; its  proof uses
a sup-norm bounded linear solution operator to the $\dibar$-equation for
$(0,1)$-forms on strongly pseudoconvex domains. 
Here we work with fiberwise injective holomorphic sprays defined 
on (not necessarily trivial) holomorphic vector bundles of class $\cA(D)$. 
To a spray over $\bar D$ we attach finitely many additional sprays over 
small open sets in $S$ whose union covers the boundary of $D$.
By improving the splitting lemma 
(Lemma \ref{splitting} below) we insure that $f(\bar D)$ 
is contained in the range of the resulting holomorphic spray; 
the inverse of this spray is the map $\Theta$ in Theorem \ref{Stein-nbds}.
We give a proof of Lemma \ref{splitting} based on the implicit 
function theorem in Banach spaces, simpler than the one in \cite{BDF1}
where iteration was used.
 
For the general theory of Stein manifolds we refer to \cite{GR} and \cite{Hor};
for real analysis in infinite dimensions see \cite{Palais1}, 
and for complex analysis in infinite dimensions see \cite{Dineen} and \cite{Lempert0}.

%
%
%  Section 2
%
%
\section{Complex manifolds of holomorphic maps}
Let $D$ be a relatively compact domain with piecewise $\cC^1$ 
boundary in $\C^n$. We consider the following function spaces on $D$:
\begin{itemize}
\item[(i)] For $k\in\Z_+$ and $0\le \alpha <1$, $\cA^{k,\alpha}(D)$
is the Banach space of all functions $\bar D\to \C$ in the H\"older class 
$\cC^{k,\alpha}(\bar D)$ which are holomorphic in $D$.
When $\alpha=0$ we shall write $\cA^{k,0}=\cA^k$ and $\cA^0=\cA$.
\item[(ii)] $\cA^\infty(D) = \cap_{k=0}^\infty \cA^{k}(D)$ 
is the Fr\'echet space consisting of all $\cC^\infty$ function $\bar D\to \C$ which 
are holomorphic in $D$.
\item[(iii)] For $k\in \Z_+$ and $p\ge 1$, $L^{k,p}_\cO(D)$ is the Banach space 
(Hilbert if $p=2$) consisting of all holomorphic functions $D\to \C$ whose 
partial derivatives of order $\le k$ belong to $L^p(D)$ (with respect to
the Lebesgue measure). 
These are Sobolev spaces of holomorphic functions on $D$.
\end{itemize}

If $L(D)$ is any of the above function spaces, we denote by $L(D,\C^m)$ 
the locally convex topological vector space consisting of maps 
whose components belong to $L(D)$. In the case (iii) we shall assume
$kp> 2n$, so the Sobolev embedding theorem provides a continuous 
(even compact) inclusion map $L^{k,p}_\cO(D) \hookrightarrow \cA(D)$
(see e.g.\ Calderon \cite{Calderon}). 

Assume now that $D$ is a relatively compact domain with piecewise $\cC^1$ boundary
in an $n$-dimensional complex manifold $S$. 
Given a complex manifold $Y$ of dimension $m$ without boundary, 
one can define the mapping space $L(D,Y)$ as follows. 
(See Lempert \cite[\S 2]{Lempert} for the case when $\bar D$ is a compact 
smooth manifold and we are considering the space $\cC^k(\bar D,Y)$,
$k\in\Z_+\cup\{\infty\}$. For the smooth manifold structure on certain 
spaces of smooth maps  see 
Palais \cite{Palais1}, as well as Cantor \cite{Cantor}, 
El\u\i asson \cite{Eliasson}, Krikorian \cite{Krikorian}, 
Penot \cite{Penot}, Riddell \cite{Riddell}.)
Fix a continuous map $f\colon \bar D \to Y$. 
Choose finitely many holomorphic coordinates systems 
$\phi_j\colon U_j\to \wt U_j\subset \C^n$ on $S$,
and $\psi_j \colon W_j \to  \wt W_j \subset \C^m$ on $Y$,
such that $f(\bar D\cap U_j) \subset W_j$ for all $j$.
Also choose open subsets $V_j\Subset U_j$ such that
$\bar D\subset \cup_j V_j$ and $V_j\cap D$ has piecewise $\cC^1$ boundary
for each $j$. 
Then $f\in L(D,Y)$ precisely when for each $j$ the restriction 
$f_j$ of the map $\psi_j\circ f\circ \phi_j^{-1}$ to the set 
$\phi_j(\overline {D\cap V_j})  \Subset \wt U_j$ 
belongs to $L(\phi_j(D\cap V_j), \C^m)$;
the definition is independent of the choices of charts. 
Further, given an open neighborhood 
$\cU_j \subset L(\phi_j(D\cap V_j),\C^m)$
of $f_j$ for every $j$, the corresponding neighborhood 
of $f$ in $L(D,Y)$ consists of all maps $g\colon \bar D\to Y$
such that $g(\overline{D\cap V_j}) \subset W_j$ and the restriction
$g_j$ of $\psi_j\circ g\circ \phi_j^{-1}$ to 
$\phi_j(\overline{D\cap V_j})$ belongs to $\cU_j$ for all $j$.

\smallskip
{\em Proof of Theorem \ref{manifolds}.}
Let $L(D,Y)$ denote any one of the above spaces;
observe that it is a subset of $\cA(D,Y)$.
We need to construct holomorphic charts in $L(D,Y)$.

Given a holomorphic vector bundle $\pi\colon E \to U$ 
over an open set $U\subset S$ containing $\bar D$,
we denote by $L_h(D,E)$ the space of all section
$\bar D\to E_{\bar D}$ of $h$ over $\bar D$ which belong to $L(D,E)$.
This is a locally convex topological vector space;
Banach for $\cA^{k,\alpha}$ or $L^{k,p}_\cO$, Hilbert for $L^{k,2}_\cO$,
and Fr\'echet for $\cA^\infty$. 

Fix a map $f\in L(D,Y)$; so $f$ is continuous on $\bar D$ and holomorphic on $D$.
Theorem \ref{Stein-nbds} furnishes an open Stein neighborhood
$\Omega\subset S\times Y$ of the graph $G_f= \{(z,f(z))\colon z\in \bar D\}$ 
and a biholomorphic map $\Theta\colon \Omega\to \wt \Omega\subset E$ 
onto an open set $\wt\Omega$ in the total space of a holomorphic vector 
bundle $\pi\colon E\to U$ such that $\bar D \subset U\subset S$
and $\pi \circ \Theta \colon \Omega\to S$ is the restriction to 
$\Omega$ of the base projection $(z,y)\to z$. 
Since $\Theta$ is holomorphic in a neighborhood of $G_f$,
the map 
\[
	\bar D\ni z \to \theta(f)(z):= \Theta(z,f(z)) \in E_z
\]
is a section of the restricted bundle $E_{\bar D}\to \bar D$ which belongs
to the space $L_h(D,E)$. The graph $G_g$ of any $g\in L(D,Y)$ sufficiently near 
$f$ is also contained in $\Omega$, and the composition with $\Theta$
defines an isomorphism $g\to \theta(g) =\Theta(\cdotp,g)$ between an 
open neighborhood of $f$ in $L(D,Y)$ and an open neighborhood 
of $\theta(f)$ in $L_h(D,E)$; we take $\theta$ as a Banach (or Fr\'echet) chart on $L(D,Y)$. 
It is easily verified that the transition map between any such pair of charts 
is holomorphic (the argument given in \cite{BDF2} for $\cA(D,Y)$ 
applies in all cases; for the Sobolev classes see \cite[Theorem 9.10]{Palais2}). 
The collection of all such charts defines a holomorphic Banach 
(resp.\ Fr\'echet) manifold structure on $L(D,Y)$.

The above construction also shows that the tangent space 
to the manifold $L(D,Y)$ at a point $f\in L(D,Y)$
can be identified with $L_h(D,f^*TY)$, the space of sections of class $L(D)$
of the complex vector bundle $h\colon f^*TY\to \bar D$
(the pull-back to $\bar D$ of the tangent bundle $TY$ by the map $f$).
By the Oka-Grauert principle, homotopic maps induce isomorphic pull-back 
bundles (see Leiterer \cite{Leiterer1} and Heunemannn \cite{Heunemann1} for
the relevant `up to the boundary' version), and hence $T_f L(D,Y)$ is independent 
of a point $f$ in a connected component of $L(D,Y)$ (up to a complex 
Banach space isomorphism). Maps belonging to different connected 
components of $L(D,Y)$ may induce nonisomorphic bundles.
%, and the tangent spaces to $L(D,Y)$ at such point need not be
%  isomorphic as complex Banach spaces. 

If $D$ is contractible, or if $\dim D=1$ 
(which means that $D$ is a bordered Riemann surface) 
then every $\cA(D)$-vector bundle over $\bar D$ is trivial; 
in this case $T_f L(D,Y) \approx L(D,\C^m)$ with $m=\dim Y$
for every $f\in L(D,Y)$.
\qed
\smallskip

The above construction is essentially the same as the one of 
Palais in the $\cC^\infty$ case \cite[Ch.\ 13]{Palais1}.
With some additional work it might be possible to introduce the categorical 
(axiomatic) approach as in \cite{Palais1} and thereby extend 
the result to a wider class of function spaces. 
Observe that the above proof actually gives the 
following more general result.

\begin{theorem}
Assume that $h\colon X\to S$ is a holomorphic submersion of a complex manifold $X$
onto a Stein manifold $S$ and $D\Subset S$ is a strongly pseudoconvex domain 
with $\cC^2$ boundary in $S$. Let $L$ be any of the classes $\cA^{k,\alpha}$,
$\cA^\infty$ or $L^{k,p}_{\cO}$ $(kp>\dim_\R S)$. 
Then the space of sections
\[
	L_h(D,X)= \{f\in L(D,X)\colon h(f(z))=z,\ z\in \bar D\}
\]
is a complex Banach manifold (Fr\'echet for $L=\cA^\infty$).
\end{theorem}

\begin{remark}
The method in \cite{BDF2} can be used to obtain the same result 
in the more general case when $h\colon X\to \bar D$ is a smooth submersion
onto $\bar D$ which is holomorphic over $D$; an example is a 
smooth fiber bundle over $\bar D$ which is holomorphic over $D$.
However, for holomorphic submersions which extend holomorphically 
to a neighborhood of $\bar D$ in $S$, the above construction is simpler 
than the one in \cite{BDF2}. Indeed, we do not need a new splitting lemma 
for each of the function spaces -- we only need it for the space 
$\cA_h(D,E)$ where $h\colon E\to\bar D$ is a complex vector
bundle which is holomorphic over $D$.
\end{remark}

The complex structures introduced above enjoy the following 
functorial property; compare with \cite[Proposition 2.3]{Lempert} 
and observe that the proof given there carries over to our situation
as well.

\begin{proposition}
Let $D\Subset S$ be a strongly pseudoconvex domain in a Stein manifold $S$,
let $Y$ and $Y'$ be complex manifolds, and let $\Phi\colon S\times Y \to Y'$
be a holomorphic map. Then the induced map $\Phi_*\colon L(D,Y) \to L(D,Y')$  
defined by $\Phi_*(f)(z)=\Phi(z,f(z))$ $(z\in \bar D)$ is holomorphic.
\end{proposition}

The discussion following Proposition 2.3 in \cite{Lempert}
also shows that the functorial property described in the above
proposition characterizes the complex structures under
consideration. 

These complex structures are functorial also with respect
to maps of the source domains:  
Given strongly pseudoconvex domains $D\Subset S$, $D'\Subset S'$
in Stein manifolds $S$ resp.\ $S'$, and given a holomorphic map 
$\Psi\colon S\to S'$ satisfying $\Psi(D)\subset D'$,
the pull-back map $\Psi^*\colon L(D',Y)\to L(D,Y)$ defined by
$f\to f\circ \Psi$ is holomorphic (see \cite[Proposition 2.4]{Lempert}). 
In particular, the evaluation map
$e\colon \bar D \times L(D,Y) \to Y$, $e(z,f) = f(z)$,
is continuous, and is holomorphic in $f$ for a fixed $z\in \bar D$ 
(see the proof of Proposition 2.5 in \cite{Lempert}). 
Weaker hypothesis on $\Psi$ may suffice in individual cases.

\section{A splitting lemma}
In this section we prove a splitting lemma for fiberwise
injective holomorphic maps on holomorphic vector bundles with
continuous boundary values. Lemma \ref{splitting} below
is the main ingredient in the proof of Theorem \ref{Stein-nbds}.
Although it is a minor extension of Theorem 3.2 in \cite{BDF1} 
(which applies to trivial bundles), we give a simpler proof 
based on the implicit function theorem. 

%
%
%  Definition of a Cartan pair
%
%
\begin{definition}
\label{Cartan-pair}
{\rm \cite[Def.\ 2.3]{BDF2}}
A pair of open subsets $D_0,D_1 \Subset S$ in a Stein manifold $S$
is said to be a {\em Cartan pair} of class $\cC^\ell$ $(\ell\ge 2)$ if 
\begin{itemize}
\item[(i)]  $D_0$, $D_1$, $D=D_0\cup D_1$ and $D_{0,1}=D_0\cap D_1$ 
are strongly pseudoconvex domains with $\cC^\ell$ boundaries, and 
\item[(ii)] $\overline {D_0\backslash D_1} \cap \overline {D_1\backslash D_0}=\emptyset$ 
(the separation property). 
\end{itemize}
{\em $D_1$ is a convex bump on $D_0$} if in addition 
there is a biholomorphic map from an open neighborhood of $\bar D_1$ in $S$ onto 
an open subset of $\C^n$ $(n=\dim S)$ which maps $D_1$ and $D_{0,1}$
onto strongly convex domains in~$\C^n$.
\end{definition}

Suppose that $D\Subset S$ is a strongly pseudoconvex domain in a Stein manifold $S$
and $\pi\colon E\to \bar D$ is a continuous complex vector bundle which is
holomorphic over $D$ (an $\cA(D)$-vector bundle). 
Such $E$ embeds as an $\cA(D)$-vector subbundle $E'$ of a trivial bundle 
$\T^N=\bar D\times\C^N$ for a sufficiently large integer $N$, 
and for every such embedding there is a Whitney 
direct sum decomposition $\T^N=E'\oplus E''$ of class $\cA(D)$. 
(These facts follows from Cartan's Theorem B for $\cA(D)$-vector bundles; 
see Leiterer \cite{Leiterer1,Leiterer2} and 
Heunemann \cite{Heunemann2,Heunemann3}.) 
Fix such a decomposition and identify $E$ with $E'$.
We shall denote the variable in $\bar D$ by $z$ and the variable in 
$\C^N$ by $t$. On $\{z\}\times\C^N$ we have a unique decomposition
$t=t'\oplus t'' \in E'_z \oplus E''_z$ which is of class $\cA(D)$ 
with respect to $z\in \bar D$.
Let $|\cdotp|$ denote the standard Euclidean norm on $\C^N$,
and let $\B=\{t\in \C^N\colon |t|<1\}$.
For every $r>0$ and $z\in \bar D$ set
\[
		E_{z,r} = \{ t= t'\oplus 0''\in E_z \colon |t|<r\} = E_z \cap r\B.
\]
Given a subset $K\subset \bar D$ and $r>0$ we shall write
\begin{equation}
\label{EK}
		E_K= \cup_{z\in K} E_z,\quad 
		E_{K,r}= \cup_{z\in K} E_{z,r} = E_K \cap (K\times r\B).  
\end{equation}
Every continuous fiber-preserving map $E_{K,r}\to E_K$
is of the form $\gamma(z,t')=(z,\psi(z,t'))$ 
for $z\in K$ and $|t'|=|t'\oplus 0''| <r$;
we shall say that $\gamma$ is of class $\cA$ if it is 
holomorphic in the interior of its domain. 
Let $id(z,t')=(z,t')$ denote the identity map on $E$. Set
\[
		||\gamma-id||_{K,r}= \sup\{ |\psi(z,t')-t'|\colon z\in K,\ |t'|<r\}.
\]

\begin{lemma}
\label{splitting}
Let $D=D_0\cup D_1$ be a Cartan pair of class $\cC^2$ in a Stein manifold $S$,
and let $\pi\colon E\to \bar D$ be an $\cA(D)$-bundle. 
Set $K=\bar D_{0,1}=\bar D_0\cap \bar D_1$. 
Given numbers $0<r'<r$ and $\epsilon>0$, there is a number $\delta>0$
satisfying the following. 
For every fiber preserving map $\gamma \colon E_{K,r}\to E_K$
of class $\cA(E_{K,r})$ with $||\gamma-id||_{K,r} <\delta$ there exist 
injective fiber preserving maps $\alpha\colon E_{\bar D_0,r'}\to E_{\bar D_0}$,
$\beta\colon E_{\bar D_1,r'}\to E_{\bar D_1}$, of class $\cA$
on their respective domains, satisfying 
$||\alpha-id||_{\bar D_0,r} <\epsilon$, $||\beta-id||_{\bar D_1,r} <\epsilon$
and 
\[
		\gamma\circ\alpha = \beta \quad {\rm on\ } E_{K,r'}. 
\]
If in addition $\gamma$ preserves the zero section (i.e., $\gamma(z,0)=(z,0)$
for $z\in K$) then $\alpha$ and $\beta$ can be chosen to satisfy the same property.

If the Cartan pair $D=D_0\cup D_1$ is of class $\cC^\ell$, $\ell\ge 2$,
then for any integer $l\in \{0,1,\ldots,\ell\}$ the analogous result holds
with the bundle $E$ and the maps $\alpha$, $\beta$, $\gamma$ of class 
$\cA^l$ on their respective domains (i.e., holomorphic inside and of
class $\cC^l$ up to the boundary), with $\cC^l$ estimates.  
\end{lemma}

\begin{proof}
We begin with the special case when the bundle $E$ is trivial, 
$E=\T^N=\bar D\times \C^N$ for some $N\in\N$. Although in this case 
the result coincides with Theorem 3.2 in \cite{BDF1},
we give a new proof based on the implicit function theorem.
(It is similar to the proof of Proposition 5.2 in \cite[p.\ 141]{FP1}.)

Recall that $(\gamma(z,t)=(z,\psi(z,t))$, where $\psi\colon K\times r\B \to \C^N$ 
is close to the map $\psi_0(z,t)=t$. 
We denote by $C_r$ (resp.\ by $\Gamma_r$) the Banach space consisting of all 
continuous maps $K\times r\B \ni (z,t)\to \psi(z,t)\in\C^N$ which are holomorphic in 
the interior $D_{0,1} \times r\B$ of $K\times r\B$ and satisfy
\begin{eqnarray*}
	||\psi||_{C_r}      &=& \sup_{(z,t)\in {K\times r\B}} |\psi(z,t)| < +\infty, \\
	||\psi||_{\Gamma_r} &=& \sup_{(z,t)\in {K\times r\B}} \bigl( |\psi(z,t)| 
									           + |\di_t \psi(z,t)| \bigr) < +\infty.									           
\end{eqnarray*}
Here, $\di_t$ denotes the partial differential with respect to the
variable $t\in\C^N$, and $|\di_t \psi(z,t)|$ is the Euclidean operator norm. 

Replacing the number $r>0$ in Lemma \ref{splitting} with a slightly 
smaller number we can assume (in view of the Cauchy estimates) 
that $\psi$ belongs to $\Gamma_r$ and that 
$||\psi-\psi_0||_{\Gamma_r}$ is as small as desired,
where $\psi_0(z,t)=t$. Fix such $r$ and choose a 
number $r'$ with $0<r'<r$.
Let $A_{r'}$ (resp.\ by $B_{r'}$) denote the Banach space of 
all continuous maps $\bar D_0\times r' \B \to \C^N$ 
(resp.\ $\bar D_1\times r' \B \to \C^N$) which are holomorphic
in the interior of the respective set, endowed with the
sup-norm. By \cite[Lemma 3.4]{BDF1} there exist continuous linear operators 
$\cA\colon C_{r'} \to A_{r'}$, $\cB\colon C_{r'}\to B_{r'}$ 
satisfying 
\begin{equation}
\label{linear}
		c=\cA(c) - \cB(c),\quad c\in C_{r'}.
\end{equation}
The proof in \cite{BDF1} uses a linear solution operator for the $\dibar$-equation 
on the level of $(0,1)$-forms on $D$ satisfying sup-norm estimates,
and the variables $t$ are treated as parameters. 
 
Given  $\psi\in \Gamma_r$ sufficiently near 
$\psi_0$ and $c\in C_{r'}$ near $0$, we define
\[
		\Phi(\psi,c)(z,t) = \psi(z,t+\cA(c)(z,t)) - (t+\cB(c)(z,t)),
		\quad (z,t)\in K\times r'\B.
\]
It is easily verified that $(\psi,c)\to \Phi(\psi,c)$ 
is a $\cC^1$ (even smooth) map from an open neighborhood
of the point $(\psi_0,0)$ in the Banach space $\Gamma_r\times C_{r'}$ 
to the Banach space $C_{r'}$.
(Although we are composing maps which are only continuous up to the boundary
of their respective domains in the $z$-variable, we are inserting 
$\cA(c)$ in the second variable of $\psi$, 
and $\psi$ is holomorphic with respect to that variable
on a larger domain.) 

Since $\Phi(\psi_0,c)= \cA(c)-\cB(c)=c$ by (\ref{linear}),
the implicit function theorem shows that in a neighborhood of 
$(\psi_0,0)$ in $\Gamma_r\times C_{r'}$
we can solve the equation $\Phi(\psi,c)=0$ on $c$; that is,
there is a $\cC^1$ map $\psi \to \cC(\psi) \in C_{r'}$,
defined in an open neighborhood of $\psi_0$ in $ \Gamma_r$ and 
satisfying 
\begin{equation}
\label{identity}
	\Phi(\psi,\cC(\psi))=0, \quad \cC(\psi_0)=0.  
\end{equation}
Consider the functions 
\[
		a_\psi = t + \cA\circ \cC(\psi) \in A_{r'},\quad 
		b_\psi = t + \cB\circ \cC(\psi) \in B_{r'}. 
\]
From (\ref{identity}) and the definition of $\Phi$ we obtain
\begin{equation}
\label{split}
		\psi(z,a_\psi(z,t)) = b_\psi(z,t),\quad (z,t)\in K\times r'\B.
\end{equation}		
Setting $\alpha(z,t)=(z,a_\psi(z,t))$,
$\beta(z,t)=(z,b_\psi(z,t))$, we see that (\ref{split}) 
gives $\gamma\circ\alpha=\beta$. We also get sup-norm estimates 
on $a_\psi-\psi_0$ (resp.\ on $b_\psi- \psi_0$) on $\bar D_0\times r'\B$ 
(resp.\ on $\bar D_1\times r'\B$) in terms of $||\psi-\psi_0||_{\Gamma_r}$. 
If the latter number is sufficiently small,  
the maps $a_\psi$ and $b_\psi$ are as close as desired 
to the map $\psi_0(z,t)=t$ and hence, after we shrink $r'$ slightly 
and apply again the Cauchy estimates in the $t$-variable, 
we can assume that they are fiberwise injective holomorphic. 
This proves Lemma \ref{splitting} when $E=\bar D  \times \C^N$.

The case with $\cC^l$ boundary values is obtained in the same
way by using the appropriate Banach spaces; compare with 
Theorem 3.2 in \cite{BDF1}.

It remain to consider the general case when $E$ is an $\cA(D)$-vector subbundle 
of $\T^N=\bar D\times \C^N$. As before we identify $E$ with its image 
$E'\subset \T^N$ and write $\T^N=E'\oplus E''$ (an $\cA(D)$-decomposition). 
Let $t=t'\oplus t''\in E'_z\oplus E''_z$ denote the corresponding splitting of the
fiber variable. We associate to each self-map $\gamma(z,t') =(z,\psi'(z,t'))$ 
of $E'$ the self-map 
\[
	\wt\gamma(z,t)=(z,\psi(z,t)), \quad \psi(z,t)= \psi'(z,t')\oplus t''
\]
of $\T^N$ (we added the identity map on the second summand $E''$).
If $\psi'$ is sufficiently close to the map $(z,t')\to t'$ then 
$\psi$ is close to the map $\psi_0(z,t)= t$, and the first part of the Lemma 
(for the trivial bundle) furnishes $\C^N$-valued maps 
\begin{eqnarray*}
		a(z,t) &=& a'(z,t)\oplus a''(z,t), \quad  (z,t)\in \bar D_0\times r'\B, \\ 
		b(z,t) &=& b'(z,t)\oplus b''(z,t), \quad \, (z,t)\in \bar D_1\times r'\B
\end{eqnarray*}
satisfying $\psi(z,a(z,t))=b(z,t)$ for $(z,t)\in K\times r'\B$.
Comparing the $E'$ and the $E''$ components of this identity 
at $t''=0$ we get
\[
	\psi'(z,a'(z,t'))=b'(z,t'),\quad a''(z,t')= b''(z,t').
\]
Hence the maps $\alpha(z,t')=(z,a'(z,t'))$ and $\beta(z,t')=(z,b'(z,t'))$
satisfy the conclusion of Lemma \ref{splitting}.
\end{proof}

\section{Proof of Theorem \ref{Stein-nbds}}
Let $f\colon \bar D\to X$ be a continuous section 
of a holomorphic submersion $h\colon X\to S$  such that 
$f$ is holomorphic in $D$. Set $\Sigma=f(\bar D)$.
Recall that $VT(X)=\ker dh$ is the vertical tangent bundle of $X$.

By Proposition 4.1 in \cite{BDF2} there exists a map 
$F\colon\bar D\times r\B^N \to X$ of class $\cA(D \times r\B^N)$ 
for some $r>0$ and $N\in\N$ such that for all $z\in\bar D$ we have
\begin{itemize}
\item[(i)] $F(\{z\} \times r\B^N) \subset X_z=h^{-1}(z)$,  
\item[(ii)] $F(z,0)=f(z)$, and 
\item[(iii)] the map $\di_t|_{t=0} F(z,t) \colon \C^N\to VT_{f(z)}X$  
is surjective.
\end{itemize}

A map $F$ with these properties is called a {\em dominating spray
of class $\cA(D)$ with the central section $F_0= F(\cdotp,0)=f$.}
(We emphasize that the sprays used here are local with respect to the parameter
and should not be confused with the global sprays used in the
Oka-Grauert theory.) 

Set $E''_z =\ker \di_t|_{t=0}F(z,t) \subset \C^N$
$(z\in \bar D)$; this defines an $\cA(D)$-subbundle 
$E''\subset \T^N :=\bar D\times \C^N$. By \cite{Heunemann2} 
and \cite{Leiterer2} there exists a complementary 
$\cA(D)$-subbundle $E'\subset \T^N$ such that $\T^N=E'\oplus E''$.
By \cite{Heunemann1} (see also the Appendix in \cite{BDF1})
we can approximate $E'$ sufficiently well over $\bar D$ 
by a holomorphic vector subbundle $E\subset U\times \C^N$
over an open set $U\subset S$ containing $\bar D$
such that $\T^N=E_{\bar D}\oplus E''$.

Let $G$ denote the restriction of $F$ to 
$E_{\bar D,r} = E_{\bar D} \cap (K\times r\B)$ (\ref{EK}).
Then $\di_{t'}|_{t'=0} G(z,t') \colon E_z \to VT_{f(z)} X$
is a linear isomorphism for every $z\in\bar D$, and by decreasing 
$r>0$ we can insure that $G$ is injective holomorphic on each fiber;
such $G$ will be called a {\em fiberwise biholomorphic spray}.
Note that $E_{\bar D}$ is isomorphic to 
the bundle $VT(X)|_\Sigma$;  when $X=S\times Y$
and $\Sigma$ is the graph of an $\cA(D,Y)$-map $f\colon \bar D\to Y$,
the latter bundle is isomorphic to $f^*TY$.

\begin{lemma}
\label{approx-sprays}
There exist a number $r' \in (0,r)$, a decreasing sequence of 
open sets $O_1\supset O_2\supset\cdots$ in $S$ with 
$\cap_{s=1}^\infty O_s=\bar D$, and a sequence of fiberwise biholomorphic sprays 
$G_s\colon E_{O_s,r'} \to X$ such that $G_s$ converges
to $G$ uniformly on $E_{\bar D,r'}$ as $s\to\infty$ and 
$\Sigma\subset G_s(E_{O_s,r'})$ for $s=1,2,\ldots$.  
\end{lemma}

\smallskip
{\em Proof of Theorem \ref{Stein-nbds}.}
Assume Lemma \ref{approx-sprays} for the moment.
Let $\Omega_0\subset X$ be an open neighborhood of $f(\bar D)$.
Choose an initial fiberwise biholomorphic spray 
$G\colon E_{\bar D,r}\to X$ as above 
such that $G(E_{\bar D,r})\subset\Omega_0$
(the latter is achieved by decreasing the number $r>0$ if
necessary). Let $G_s\colon E_{O_s,r'} \to X$ $(s=1,2,\ldots)$ 
be a sequence of sprays furnished by Lemma \ref{approx-sprays}. Set 
\[
	\wt \Omega_s = E_{O_s,r'} \subset E, \quad 
	\Omega_s = G_s(E_{O_s,r'}) \subset X,
	\quad \Theta_s=G_s^{-1} \colon \Omega_s\to\wt \Omega_s. 
\]
If $s\in\N$ is chosen sufficiently large then $\Omega_s\subset \Omega_0$,
and for such $s$ the map $\Theta_s\colon \Omega_s\to\wt \Omega_s$ 
satisfies the conclusion of Theorem \ref{Stein-nbds}.
\qed

\smallskip
{\em Proof of Lemma \ref{approx-sprays}.}
We proceed as in the proof of Theorem 5.1 in \cite{BDF2},
but paying close attention to the ranges of the 
approximating sprays in order to insure that each of them
contains the graph $\Sigma=f(\bar D)$. 

Since $h\colon X\to S$ is a holomorphic submersion,  there exist
for each point $x_0\in X$ open neighborhoods $x_0\in W\subset X$,
$h(x_0)\in V\subset S$, and biholomorphic maps 
$\phi \colon V\to \B^n \subset \C^n$,
$\Phi \colon W\to \B^n\times \B^m \subset \C^n\times\C^m$,
such that $\phi(h(x))=pr_1(\Phi(x))$ for every $x\in W$.
Note that 
\[
	\Phi(x)= \bigl(\phi(h(x)),\phi'(x)\bigr) \in \B^n\times \B^m, \quad x\in W,
\]
where $\phi'=pr_2\circ\Phi$. We call such $(W,V,\Phi)$ 
a {\em special coordinate chart} on $X$.

Recall that $G\colon E_{\bar D,r} \to X$ is fiberwise biholomorphic spray 
over $\bar D$ with the central section $f$.
Nara\-sim\-han's lemma on local convexification of a strongly pseudoconvex 
hypersurface gives finitely many special coordinate charts 
$(W_j,V_j,\Phi_j)$ on $X$, with $\Phi_j=(\phi_j\circ h,\phi'_j)$,
such that $bD\subset \cup_{j=1}^{j_0} V_j$ 
and the following hold for all $j=1,\ldots,j_0$
(for (ii) and (iii) we may have to decrease $r>0$):
\begin{itemize}
\item[(i)]   $\phi_j(bD\cap V_j)$ is a strongly convex hypersurface 
in the ball $\B^n$, 
\item[(ii)]  the spray $G$ maps $E_{\bar D\cap V_j,r}$ into $W_j$, and 
\item[(iii)] $\phi'_j\circ G(E_{\bar D\cap V_j,r}) \Subset \B^m$.
\end{itemize}

Fix an $r>0$ such that the above properties hold and 
choose a number $r'\in (0,r)$. Also choose a number $c\in (0,1)$ sufficiently close to $1$ 
such that the open sets 
$U_j= \phi_j^{-1}(c\B^n) \Subset V_j$ $(j=1,\ldots,j_0)$ still cover $bD$.

By a finite induction we shall find strongly pseudoconvex domains 
$D=D_0\subset D_1\subset\cdots\subset D_{j_0} \Subset U$,
numbers $r=r_0>r_1>\cdots > r_{j_0}=r'$ and fiberwise biholomorphic sprays 
$G_k \colon E_{\bar D_k, r_k} \to X$ of class $\cA$
$(k=0,1,\ldots,j_0)$, with $G_0=G$, 
such that for every $k\in \{1,\ldots, j_0\}$ the restriction 
of $G_k$ to $E_{\bar D_{k-1},r_k}$ will approximate $G_{k-1}$ 
in the sup-norm topology. The domain $D_k$ will be chosen such that 
\[
	D_{k-1}\subset D_k \subset D_{k-1}\cup V_k, \quad bD_{k-1} \cap U_k\subset D_k
\]
for $k=1,\ldots,j_0$; that is, we enlarge (bump out) $D_{k-1}$ 
inside $V_k$ so that the part of $bD_{k-1}$ which lies in the smaller set $U_k$ 
is contained in the next domain $D_k$. As the $U_j$'s cover $bD$,
the final domain $D_{j_0}$ will contain $\bar D$ in its interior, 
and the spray $\wt G= G_{j_0} \colon E_{\bar D_{j_0},r'} \to X$ will 
approximate $G$ as close as desired on $E_{\bar D,r'}$;
in particular, we shall arrange that 
$\Sigma=f(\bar D)$ is contained in $\wt G(E_{\bar D_{j_0},r'})$.
To keep the induction going we will also insure at every step
that the properties (ii) and (iii) above remain valid with 
$(D,G)$ replaced by $(D_k,G_k)$ for all $k=1,\ldots,j_0$.
The restriction of $\wt G$ to the interior $E_{D_{j_0},r'}$
can be taken as one of the sprays in the conclusion of the lemma.

\begin{figure}[ht]
\psset{unit=0.6cm, linewidth=0.7pt}  
\begin{pspicture}(-8,-6)(8,6)

\pscircle[linecolor=OrangeRed,linewidth=1pt,fillstyle=none](0,0){6}
\pscircle[linecolor=Red,linewidth=1pt,linestyle=dashed](0,0){3.2}

\pscurve(-3.35,5.4)(1.5,0)(-3.35,-5.4)                               % The boundary of $D_1$
\psecurve[linecolor=DarkBlue](-4,5)(-1.4,4)(2.4,0)(-1.4,-4)(-4,-5)   % The boundary of $D'_1$
\psecurve[linecolor=Purple](-4,5)(-1.4,4)(0.6,0)(-1.4,-4)(-4,-5)     % The boundary of $D$
\psarc(0,0){5.4}{132}{228}                                           % arc - boundary of $D'_1\cap D$                         
\psecurve(-4,3.4)(-3.6,4)(-2.5,4.5)(-1.9,4.4)(-1.4,4)                % top connecting curve
\psecurve(-4,-3.4)(-3.6,-4)(-2.5,-4.5)(-1.9,-4.4)(-1.4,-4)           % borrom connecting curve 

\rput(0,-5){$bD$}
\psline[linewidth=0.2pt]{<-}(-2.63,-5)(-0.5,-5)

\rput(-1.2,0){$bD\cap U_1$}
\psline[linewidth=0.2pt]{->}(-0.1,0)(0.6,0)

\rput(4.6,0){$bD'_1$}
\psline[linewidth=0.2pt]{->}(4,0)(2.45,0)

\rput(3.7,-2.8){$bD_1$}
\psline[linewidth=0.2pt]{->}(3.1,-2.6)(1.25,-1)

\rput(3.2,3){$U_1$}
\psline[linewidth=0.2pt]{->}(2.8,2.8)(2.3,2.3)

\rput(0.7,5){$V_1$}

\rput(-7.8,0){$D\cap D'_1$}
\psline[linewidth=0.2pt]{->}(-6.6,0)(-4.8,0)

\end{pspicture}
\caption{The domains $D'_1$ and $D_1$}
\label{Fig1}
\end{figure}
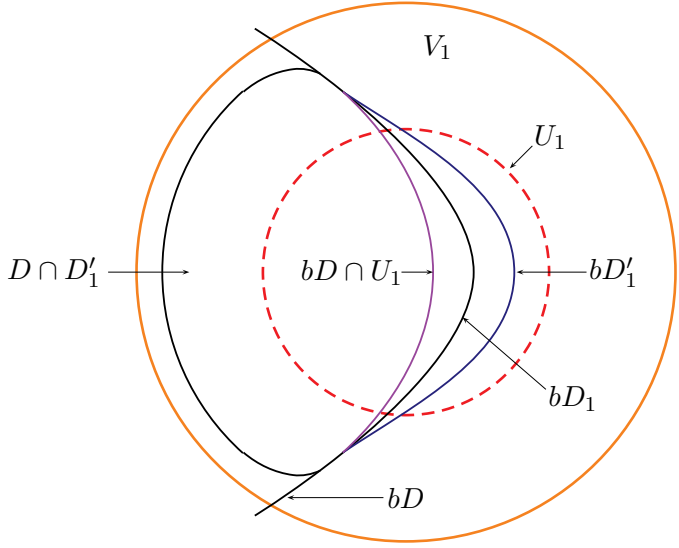

The geometric scheme is as in the proof of Theorem 5.1 in \cite{BDF2}.
Since all steps are of the same kind, it suffices to 
explain how to get the pair $(D_1,G_1)$ from $(D,G)=(D_0,G_0)$. 
We begin by finding a domain $D'_1\subset S$ with $\cC^2$ boundary 
which is a convex bump on $D=D_0$ (Definition \ref{Cartan-pair}) and 
such that $\overline  U_1 \cap \bar D \subset \bar D'_1 \subset V_1$.
To do this, we shall first find a set $\wt D'_1 \subset \B^n$ 
with suitable properties and then take $D'_1=\phi_1^{-1}(\wt D'_1)$.
Choose a smooth function $\chi\ge 0$ with compact support on $\B^n$  
such that $\chi=1$ on $c\B^n$. Recall that $U_1=\phi_1^{-1}(c\B^n)$.
Let $\tau\colon \B^n \to\R$ be a strongly convex defining
function for the domain $\phi_1(D \cap V_1) \subset \B^n$.
Choose a number $c'\in (c,1)$ close to $1$ such that the hypersurface
$\phi_1(bD\cap V_1)=\{\tau=0\}$ intersects the sphere 
$\{\zeta\in\C^n\colon |\zeta|=c'\}$ transversely. 
If $\delta >0$ is sufficiently small then the set  
\[
	\{\zeta \in \C^n \colon |\zeta| < c',\ \tau(\zeta) < \delta \chi(\zeta)  \}
\]
could serve our purpose, except that it is not
smooth along the intersection of the (convex) hypersurfaces $\{|\zeta| = c'\}$
and $\{\tau = \delta \chi \}$. By rounding off the corners
we get a strongly convex set $\wt D'_1 \subset \B^n$ 
such that $D'_1 = \phi_1^{-1}(\wt D'_1)\subset V_1$ satisfies the
desired properties. (See Fig.\ \ref{Fig1} which is 
taken from \cite{BDF2}.)  

%
%
%  Figure 1: The domains $D'_1$ and $D_1$
%
%

Choose numbers $r_1, r'_1, r''_1$ with $r' < r_1 < r'_1 < r''_1 < r$. 
By using the special coordinate chart $(W_1,V_1,\Phi_1)$ 
we find an open set $V'_1  \subset V_1$ containing 
$\bar D\cap V_1$ ($V'_1$ will depend on the choice of $G'$)
and a fiberwise biholomorphic spray $G' \colon E_{\bar V'_1, r''_1} \to X$
with range in $W_1$ whose restriction to $E_{\bar D\cap \bar V_1, r''_1}$
approximates the initial spray $G=G_0$ as close as desired in the uniform topology.
If the approximation is sufficiently close on $E_{\bar V'_1, r''_1}$,
there exists a (unique) fiberwise biholomorphic map 
$\gamma\colon E_{\bar D\cap \bar V_1, r'_1} \to E$ of class $\cA$
which is close to the identity map and satisfies the equation
\[
		G(z,t) = G'(\gamma(z,t)) = G'(z,\psi(z,t)), \quad 
		z\in \bar D\cap \bar V_1,\ t\in E_{z,r'_1}.
\]
Applying Lemma \ref{splitting} on the  Cartan pair 
$(D,D'_1)$ we obtain $\gamma\circ\alpha =\beta$ 
on $E_{\bar D\cap \bar V_1, r_1}$, where $\alpha\colon E_{\bar D,r_1} \to E$ 
and $\beta \colon E_{\bar D'_1,r_1} \to E$ are injective holomorphic maps 
which are close to the identity on their respective domains. 
It follows that the fiberwise biholomorphic sprays
\[
	G\circ \alpha \colon E_{\bar D,r_1} \to X, \quad 
	G'\circ \beta \colon E_{\bar D'_1\cap \bar V'_1, r_1}\to X 
\]
agree on the intersection of their domains, and hence they
define a fiberwise biholomorphic spray 
$G_1\colon E_{\bar D\cup (\bar D'_1\cap \bar V'_1), r_1}\to X$ 
of class $\cA$. By construction $G_1$ approximates $G$ uniformly on 
$E_{\bar D,r_1}$ as close as desired.

It remains to restrict $G_1$ to a suitably chosen strongly pseudoconvex
domain $D_1\Subset S$ contained in $D\cup (D'_1\cap V'_1)$ and satisfying
the other required properties. We choose $D_1$ such that it agrees with 
$D$ outside of $V_1$, while 
\[
	D\cap V_1 =\phi_1^{-1}(\{ \zeta \in \B^n \colon \tau(\zeta) < \e \chi(\zeta)\})
\]
for a small $\e>0$ (Fig.\ \ref{Fig1}). By choosing $\e$ sufficiently small 
(depending on $G_1$) we can insure that properties (i)--(iii) are satisfied by
the pair $(D_1,G_1)$. 

Applying the same procedure to $(D_1,G_1)$ and the special coordinate chart 
$(W_2,V_2,\Phi_2)$  we get the next pair $(D_2,G_2)$.  
After $j_0$ steps we find a domain $D_{j_0}\subset S$ containing 
$\bar D$ and a fiberwise biholomorphic spray 
$\wt G=G_{j_0} \colon E_{\bar D_{j_0},r'} \to X$ which
approximates $G$ as close as desired uniformly on $E_{\bar D,r'}$.
If the approximation is close enough then the range of $\wt G$ 
contains $\Sigma$.

The sequence $G_s$ in Lemma \ref{approx-sprays} 
is chosen to consist of sprays $\wt G$ obtained as above, 
approximating $G$ ever more closely on $E_{\bar D,r'}$.  
\qed

\section{Examples and problems}
Let $\D$ denote the open unit disc in the complex plane $\C$.
Our first example shows that Corollary \ref{cor1} fails in general 
for maps which are discontinuous at the boundary.

\begin{example}
\em There exists a bounded holomorphic function on $\D$ such that the
closure of its graph does not have a Stein neighborhood basis in $\C^2$.

\rm  Indeed, let $f\in H^\infty(\D)$ be a bounded holomorphic function
on the unit disc such that $\sup_{z\in\D} |f(z)|=1$ and the cluster set of $f$
at every boundary point $e^{i\theta} \in b\D$ equals $\overline \D$.
(Such functions are easily found by using interpolation theorems for 
$H^\infty(\D)$, see \cite[Ch.\ VII]{Garnett}.)
The closure $K$ of the graph of $f$ in $\C^2$ is the union of the graph 
with all vertical discs $\{e^{i\theta}\}\times\bar \D$, $\theta\in\R$.
By the classical argument of Hartogs \cite{Hartogs} any open Stein neighborhood 
of $K$ in $\C^2$ also contains the unit bidisc. 
\end{example}

\begin{problem}
Characterize the bounded holomorphic functions on the disc $\D$ 
for which the closure of the graph admits a basis of open Stein
neighborhoods in $\C^2$.
\end{problem}

The next two examples illustrate that Corollary \ref{cor1} fails
in general for {\em images} (as opposed to graphs) of maps.
See however \cite[Theorem 2.1]{BDF1}.

%
%
%   Stout's example
%
%
\begin{example} This example was communicated to me by E.\ L.\ Stout
(private communication, September 19, 2006):

\textit{For every $N>1$ there is an $\cA(\D)$-map $\overline \D \to \C^N$ 
whose image has no Stein neighborhood basis.}

Let $E$ be a Cantor set of length zero contained in $b\D$;
then $E$ is a peak-interpolation set for the disc algebra $\cA(\D)$. 
Let $\B$ denote the unit ball in $\C^N$ for some $N>1$.
Choose a continuous map $f$ from $E$ onto the sphere $b\B$.
There is a map $F=(F_1,\ldots,F_N)\colon\bar \D \to\C^N$,
satisfying $F_j\in\cA(\D)$ for each $j$, such that $F|_E=f$
and $|F(z)|^2=\sum_{j=1}^N |F_j(z)|^2 <1$ for all $z\in\overline \D\bs E$ 
(Globevnik \cite{Globevnik} and Stout \cite{Stout}).
Now let $D'$ be a domain in $\D$ obtained by moving each of the 
open arcs of $b\D \bs E$ in just a little, leaving the end 
points fixed; so $D'$ is conformally a disc and $bD'\cap b\D = E$.  
Then $F(\bar D')$ is a compact set consisting of the sphere $b\B$ 
together with a proper subset of $\B$, and hence it has no 
Stein neighborhood basis (any Stein neighborhood also contains
the ball $\B$). 
It is necessary to pass to a smaller domain $D'\subset\D$ 
because $F$ might take $\overline \D$ onto the  ball 
$\overline{\,\B}$  which has a basis of Stein neighborhoods.
\end{example}

\begin{example}
This example is a minor modification of the one which was communicated 
to me by J.-P.\ Rosay on April 6, 2004:

\textit{There exists a smooth ($\cC^\infty$) injective map $\Phi$ from the closed unit ball 
$\overline {\,\B}$ in $\C^5$ into $\C^{8}$, 
that is a holomorphic embedding of the open unit ball, such that 
$\Phi(\overline {\,\B})$ has no basis of Stein neighborhoods.}

\smallskip
We proceed as follows. The set 
\[
	M=\{ (z_1,\ldots ,z_5)\in \overline {\,\B} \colon z_1z_2\cdots z_5 = {\sqrt 5}^{\,-5} \}
\]
is a real four dimensional submanifold of the boundary of $\B$
which is complex tangential to the sphere $b\B$ at each point.
By Chaumat and Chollet (\cite{CC1}, \cite{CC2}) every compact subset of $M$
is a peak interpolation set for $\cA^\infty(\B)$, 
the Fr\'echet algebra of functions holomorphic on the ball and smooth 
up to the boundary. Let $H$ be a closed Hartogs figure in $\C^2$:
\[
	 (\zeta_1,\zeta_2)\in H \iff
   (|\zeta_1|\leq 1~,~|\zeta_2|\leq \frac{1}{2})~{\rm or}~
( \frac{1}{2} \leq |\zeta_1|\leq 1~,~|\zeta_2|\leq 1).
\]
Let $H_0$ be a diffeomorphic copy of $H$ in $M$ (such exists by dimension reasons).
Consider a smooth map $\varphi\colon \overline{\,\B} \to \C^2$ that is holomorphic on 
the open ball and whose restriction to $H_0$ is a diffeomeorphism from $H_0$ onto $H$. 
Also choose a function $h\in \cA^\infty(\B)$ that is zero on $H_0$ and that vanishes 
nowhere else on the closed unit ball. (Such $\phi$ and $h$ are obtained by
appealing to \cite{CC1}, \cite{CC2}.)
Define a map $\Phi\colon\overline{\,\B} \to \C^8$ by 
\[
	\Phi(z)= \Phi (z_1,\ldots ,z_5) = (\phi (z), h(z),z_1h(z),\ldots, z_5h(z)).
\]
It is easy to check that $\Phi$ is injective on $\overline{\,\B}$,
it is of maximum rank (immersion) at every point of $\overline{\,\B} \bs H_0$, 
it maps $H_0$ onto $H \subset \C^2\times \{ 0\}^{6}$,
and $\Phi(z)\notin \C^2\times \{0\}^{6}$ if $z\in \overline {\,\B} \setminus H_0$.

Since $\Phi (\overline {\,\B})$ contains the Hartogs figure $H\times\{0\}^{6}$,
any open Stein neighborhood of it will also contain the unit bidisc in 
$\C^2 \times \{ 0\}^{6}$; as this bidisc  is not included in $\Phi (\overline {\,\B})$,
the latter set has no basis of Stein neighborhoods. 
\end{example}

\begin{problem}
Let $K$ be a compact set with a Stein neighborhood basis 
in a complex manifold $S$. Assume that $f\colon K\to Y$ is a continuous map 
to a complex manifold $Y$ which is a uniform limit on $K$ 
of a sequence of holomorphic maps $f_j\colon V_j\to Y$ 
defined in open neighborhoods of $K$.

\textit{Does the graph $G_f=\{(z,f(z))\colon z\in K\}$
admit a basis of open Stein neighboorhoods 
in $S\times Y$ ?}

The answer is easily seen to be affirmative if $Y=\C^N$ 
and, by the embedding theorem, also if $Y$ is 
a Stein manifold. When $K$ is the closure of a strongly 
pseudoconvex domain the answer is given by Corollary \ref{cor1}.
Another case of interest is  
the closure of a {\em weakly pseudoconvex domain} $D$ 
such that $K=\bar D$ admits a Stein neighborhood basis; 
the latter condition is necessary as is shown by the 
worm domain of Diederich and Forn\ae ss \cite{DF}.
\end{problem}

%
%
%  THANKS, THANKS AND THANKS
%
%
%

\medskip
\textit{Acknowledgements.} 
The methods used in this paper are based on my joint works 
\cite{BDF1,BDF2} with Barbara Drinovec-Drnov\v sek whom I thank
for her indirect contribution. I also wish to thank Debraj
Chakrabarti and Laszlo Lempert for helpful remarks, 
Jean-Pierre Rosay and Edgar Lee Stout
for kindly letting me include their examples in \S 5, and 
the referee for kindly explaining me the  construction
of Stein neighborhoods in Theorem \ref{Stein-nbds}
by using plurisubharmonic functions.

\bibliographystyle{amsplain}

\end{document}